%November 14, 2001
%This paper is written in AMSTeX
\documentstyle{amsppt}
\magnification=1200
\hoffset=-0.5pc
\vsize=57.2truepc
\hsize=38truepc
\nologo
\spaceskip=.5em plus.25em minus.20em
\define\BK{Barannikov-Kontsevich\ }

\define\Gammat{\tau}
\define\Sigm{\Cal S}
\define\fra{\frak}
\define\fiel{\bold}

\define\partia{\Delta}

\define\Bobb{\Bbb}

\define\barakont{1}
\define\batviltw{2}
\define\ebrown{3}
\define\cartanon{4}
\define\cartantw{5}
\define\cartanse{6}
\define\cartanei{7}
\define\chenone{8}
\define\chenfou{9}
\define\degrmosu{10}
\define\eilmactw{11}
\define\geschthr{12}
\define\getzltwo{13}
\define\gugenhtw{14}
\define\gugenhth{15}
\define\gugenlam{16}
\define\gulasta{17}
\define\gulstatw{18}
\define\gugenmun{19}
\define\hiltoone{20}
\define\homotype{21}
\define\perturba{22}
\define\cohomolo{23}
\define\modpcoho{24}
\define\intecoho{25}
\define\extensta{26}
\define\bv{27}
\define\lrbata{28}
\define\twilled{29}
\define\banach{30}
\define\berikas{31}
\define\huebkade{32}
\define\kodnispe{33}
\define\kosmathr{34}
\define\koszulon{35}
\define\ladastas{36}
\define\liazutwo{37}
\define\maninbtw{38}
\define\maninfiv{39}
\define\massey{40}
\define\moorefiv{41}
\define\munkholm{42}
\define\quilltwo{43}
\define\sanebfou{44}
\define\schlstas{45}
\define\schlsttw{46}
\define\uehamass{47}
\define\xuone{48}
\define\Nsddata#1#2#3#4#5{
 ( #4
\text{
\vbox 
to 1.15 pc
{
  \hbox {$@>{\,\,\, #3\,\,\,}>>$}
  \vskip-1.2pc
  \hbox {$@<<{\,\,\, #2\,\,\,}<$}
                }
      }
 #1,#5 )
}

\topmatter
\title
Formal solution of the master equation\\
via HPT and deformation theory
\endtitle
\author
Johannes Huebschmann and Jim Stasheff
\endauthor

\abstract
{We construct a solution of the master equation by means of standard tools from
homological perturbation theory
under just the  hypothesis
that the ground field be of  characteristic zero,
thereby avoiding
the formality assumption of the relevant dg Lie algebra.
To this end, we endow 
the homology
$\roman H(\fra g)$
of any differential
graded Lie algebra $\fra g$
over a field of characteristic zero
with  an sh-Lie structure such that
$\fra g$ and $\roman H(\fra g)$
are sh-equivalent.
We discuss our solution of the master equation in the context
of deformation theory.
Given the extra structure appropriate to the extended moduli space of
complex structures on a Calabi-Yau manifold, the known solutions result
as a special case.}
\endabstract
\affil
{USTL, UFR de Math\'ematiques,
\hfill\linebreak
Math-UNC}
\endaffil

\address{\noindent
USTL, UFR de Math\'ematiques, F-59 655 Villeneuve d'Ascq C\'edex,
France
\newline\noindent
Johannes.Huebschmann\@agat.univ-lille1.fr
\newline\noindent
Math-UNC, Chapel Hill NC 27599-3250, USA,
jds\@math.unc.edu}
\endaddress
\subjclass{
13D10
14J32
16S80  
16W30 
17B55 
17B56 
17B65
17B66
17B70  
17B81
18G10 
32G05
55P62
55R15
81T70}
\endsubjclass
\keywords
{Twisting cochain, deformation equation, master equation,
classification of rational homotopy types, formal power
series connection, differential graded algebra, differential graded
Lie algebra, sh-Lie algebra,
homological perturbation theory, dG algebra,
dBV algebra, Frobenius structure}
\endkeywords
\endtopmatter
\document
\rightheadtext{Solution of the master equation and HPT}

\beginsection Introduction

The name \lq\lq Master Equation\rq\rq\ 
derives from the physics literature, especially that of the
Batalin-Vilkovisky approach to Lagrangians with symmetries
\cite\batviltw,
but the master equation has many precursors in mathematics,
among which the best known is perhaps the Maurer-Cartan
equation. The master equation makes sense as an equation on
elements of a differential graded algebra, associative or Lie
or suitable higher homotopy generalizations. As an equation 
on elements of a differential graded Lie algebra $\Cal L$, for an element
$\Gammat$ of $\Cal L$, the master equation has the form
$$
D \Gammat = \frac 12 [\Gammat,\Gammat].
$$
In the literature, it is also customary to take the equation
$D \Gammat + \frac 12 [\Gammat,\Gammat]=0$
as (Lie algebra) master equation, the choice of sign being a matter
of convention.
As an equation on elements of a differential graded associative algebra
$\Cal A$, the master equation has the form
$D \Gammat = \Gammat\Gammat$.
By the methods of {\it homological perturbation theory\/}, we will construct 
formal solutions of the master equation in the following general context.
Consider a 
differential graded Lie algebra $\fra g$,
and denote by
$\Sigm^{\roman c}[s\fra g]$
the differential graded symmetric coalgebra
on 
the suspension $s\fra g$ of $\fra g$;
then
$\roman{Hom}(\Sigm^{\roman c}[s\fra g], \fra g)$
inherits a differential graded Lie algebra structure and, 
in view of the defining property of a Lie algebra twisting cochain
(see Section 1 below for details),
the {\it universal twisting cochain\/}
$\tau_{\fra g}\in \roman{Hom}(\Sigm^{\roman c}[s\fra g], \fra g)$
satisfies the equation
$D\tau_{\fra g} = \frac 12 [\tau_{\fra g},\tau_{\fra g}]$
which is plainly just the master equation.
By means of our main result, Theorem 2.7 below,
up to s(trong) h(omotopy) equivalence,
we can then replace 
the differential graded coalgebra $\Sigm^{\roman c}[s\fra g]$
with a differential
graded coalgebra of the kind
$\Sigm^{\roman c}_{\Cal D}[s\roman H( \fra g)]$
where 
$\roman H( \fra g)$ is the homology of $\fra g$ and
$\Cal D$  a coalgebra differential
on
$\Sigm^{\roman c}[s\roman H( \fra g)]$
turning the latter into a coaugmented differential graded 
coalgebra---the differential $\Cal D$ defines an
sh-Lie structure on $\roman H( \fra g)$, 
see Section 2 below---in such a way that 
the following hold:
The differential graded Lie algebra
$\Cal L=\roman{Hom}(\Sigm^{\roman c}_{\Cal D}[s\roman H( \fra g)], \fra g)$
is sh-equivalent to
$\roman{Hom}(\Sigm^{\roman c}[s\fra g], \fra g)$,
and the 
data give rise to
a twisting cochain
$\tau \in\roman{Hom}(\Sigm^{\roman c}_{\Cal D}[s\roman H( \fra g)], \fra g)$
which, in a suitable sense, is equivalent to
$\tau_{\fra g}$.
The twisting cochain $\tau$ is our most general
solution of the master equation.
Our 
approach in terms of homological perturbation theory will
yield
explicit recursive formulas for
$\Cal D$ and $\tau$,
once
a choice of contraction 
$\Nsddata {\fra g} {\nabla}{\pi}{\roman H(\fra g)}h$
of chain 
complexes 
has been made. This is always possible over a field.
(See Section 2 and especially (2.7.3) and (2.8.2) below for details.)
Our Theorem 2.7
{\it also\/} establishes the fact (which has been known to both of us for 
some time)
that, roughly speaking, sh-Lie structures are preserved
under strong deformation retractions; in particular,
given a differential graded Lie algebra $\fra g$,
its homology $\roman H(\fra g)$ admits an sh-Lie structure
such that $\fra g$ and
$\roman H(\fra g)$ are sh-equivalent. 

A situation of particular interest recently is that in which
$\fra g$ is in fact a differential Gerstenhaber algebra.
In the case of the moduli space of complex structures
on a complex manifold $M$, the relevant differential
Gerstenhaber algebra contains the ordinary
Kodaira-Spencer (differential graded Lie)
algebra which, in turn, is well known to control
deformations of the complex structure.
When $M$ is a Calabi-Yau manifold,
\BK \cite\barakont\ 
construct a formal solution
of the {\it master equation\/}.
They use this solution to construct a formal Frobenius
manifold structure on the extended moduli space
of complex structures on $M$.
The \BK formal solution of the master equation
results from our construction
as a special case, see Theorem 3.8 below.
Thus our  homological perturbation theory
construction of the twisting cochain $\tau$
extends
the known constructions of
solutions of the master equation
to much less restrictive hypotheses
(indeed, over a field of characteristic zero,
we do not need a hypothesis at all)
and, in particular,  allows for generalization 
to the case in which
the underlying manifold does not admit a K\"ahler structure.

We have already observed that the master equation
recovers the defining property of
a twisting cochain.
Twisting cochains have a history in topology and differential
homological algebra \cite{\ebrown,\, \cartanse,\,\moorefiv}.
Certain twisting cochains were studied by Chen
\cite{\chenone,~\chenfou} 
under the name
{\it power series connection\/}.
Berikashvili and his students 
have studied the classification of fibrations with fixed base
and fiber using
Berikashvili's functor
$D$ which is defined
in terms of homotopy classes of twisting cochains
(or twisting elements);
see \cite \sanebfou\ (Section 2)
or \cite\berikas\ 
for a recent account of Berikashvili's functor
$D$ and references to work of
Berikashvili, Kadeishvili, and Saneblidze in \cite{\berikas,~\sanebfou}.
Twisting cochains play a crucial role in
the deformation theory for rational homotopy types
and rational fibrations 
and for the corresponding classification theory
developed by Schlessinger and Stasheff
in \cite\schlstas, as well as for the construction
of small resolutions for doing calculations
in group cohomology
\cite{\perturba--\intecoho}.

That twisting cochains
may be constructed by means of  homological perturbation
theory
has been known for some time,
see e.~g.
\cite\gugenmun, 
\cite\perturba\ (2.11),
\cite\munkholm\ (2.2).
More historical comments about homological perturbation
theory may be found e.~g. in Section 1
(p. 248) and Section 2 (p. 261) of \cite\huebkade,
which has one of the strongest results
in relation to compatibility
with other 
(such as algebra or coalgebra) structure,
since it was perhaps first recognized in \cite\perturba.
Homological perturbation theory
constructions enabled the first-named
author to carry out complete numerical calculations
in group cohomology \cite\cohomolo,
\cite\modpcoho,
\cite\intecoho\ 
which cannot be done
by other methods.
This was an instance where 
homological perturbation theory led
to results which are independent of this theory.
The solution $\tau$ of the master equation
and its offspring
constructed in the present paper
constitute
as well a result which is phrased
independently of homological perturbation theory.
The significance of our more general construction for
deformation theory 
is not yet fully understood.
It is interesting to observe, though,
that the 
inductive construction
of the formal power series
$\varphi(t)$
in \S 3 of \cite\kodnispe\ 
may be seen as a homological perturbation theory
construction.
Some applications 
related with 
\cite\schlstas\ 
and generalizations will be given in Section 4.

%\beginsection 
\medskip\noindent{\bf 1. Master equation and twisting cochains}
\smallskip\noindent
The ground ring is 
assumed to contain the rationals as a subring and will be
written $R$.
We will take {\it chain complex\/}
or {\it dg module\/}
to mean a differential graded $R$-module.
A chain complex will not necessarily be concentrated
in non-negative (or non-positive) degrees.
In the present section and in Section 2,
the differential $d$ on a chain complex $M$ will be considered
as lowering degree by $1$
and referred to as a {\it homology\/} differential
(if need be);
a {\it cohomology differential\/} 
$\delta \colon M^j \to M^{j+1}$
may always be written in this way as
$d \colon M_{-j} \to M_{-j-1}$
where
$M_*= M^{-*}$ etc.
and a homology differential may accordingly be written
as a cohomology differential.
We will come back to this in Section 3 below.

Write $s$ for the {\it suspension\/} operator as usual
and accordingly $s^{-1}$ for the
{\it desuspension\/} operator.
Thus $(sM)_j = M_{j-1}$, etc.
The algebraic geometer's notation 
for the suspended (or desuspended, depending on terminology) object
is
$\Cal A[1]$ or $\Cal A[-1]$
(depending on the convention of grading).
We avoid this notation.

The coaugmentation map of a coaugmented differential graded
coalgebra $C$ will be denoted by $\eta \colon R \to C$.
The graded tensor algebra,
graded tensor coalgebra,
and graded symmetric coalgebra functors
will be written 
respectively $\roman T$,
$\roman T^{\roman c}$, $\Sigm^{\roman c}$;
see what is said below for details.
Given two chain complexes
$X$ and $Y$, 
recall that 
$\roman{Hom}(X,Y)$
inherits the structure of a chain complex
by the operator $D$  
defined by
$$
D \phi = d \phi -(-1)^{|\phi|} \phi d
$$
where $\phi$ is a homogeneous homomorphism
from $X$ to $Y$ and where $|\phi|$ refers to its degree.
\smallskip\noindent
{\smc Definition 1.1.}
A {\it dg Lie algebra\/} $\fra g$ is a dg module
$\fra g = \{\fra g_n\}$
with a differential $d\colon \fra g_n \to \fra g_{n-1}$ and
a graded degree zero bracket
$[\ , \ ]\colon \fra g_p \otimes \fra g_q \to \fra g_{p+q}$
which is a chain map and satisfies
the graded Jacobi identity, that is, for homogeneous
$X \in \fra g$, the operation
$[X,\ ]$ is a graded derivation for the bracket.

\smallskip\noindent
{\smc Definition 1.2.}
For a filtered dg module $X$, 
a {\it perturbation\/} 
of the differential $d$ of $X$ 
is a (homogeneous) morphism
$\partial$ of the same degree as
$d$
such that $\partial$ lowers the filtration and
$(d + \partial)^2 = 0$
or, equivalently,
$$
[d,\partial] + \partial \partial = 0.
\tag1.2.1
$$

Thus, when $\partial$ is a perturbation on $X$,
the sum $d + \partial$,
referred to as the {\it perturbed differential\/},
endows $X$
with a new differential.
When $X$ has a graded coalgebra structure
such that $(X,d)$ is a differential graded coalgebra,
and when the
{\it perturbed differential\/} $d + \partial$
is compatible with the graded coalgebra structure,
we refer to $\partial$ as a {\it coalgebra perturbation\/};
similarly, we can also talk about
an {\it algebra perturbation\/}.

Consider a more \lq sophisticated\rq\ 
description of dg Lie algebra:
Let $M$ be a chain complex. Write
$$
\roman T^{\roman c}[M] = (\oplus_{k \geq 0} M^{\otimes k}, \Delta,d)
$$
for the {\it coaugmented differential graded tensor coalgebra\/}
on $M$; here $\Delta$ and $d$ refer to the diagonal
and differential respectively, $M^{\otimes 0} = R$,
the coaugmentation map $\eta\colon R \to \roman T^{\roman c}[M]$
is the canonical one, and, for every $n \geq 0$,
$$
\Delta(m_1 \otimes \dots \otimes m_n) 
= \sum_{p=0}^n
(m_1 \otimes \dots \otimes m_p)\otimes(m_{p+1} \otimes \dots \otimes m_n),
$$
where $m_1,\dots,m_n \in M$.
The differential graded
symmetric coalgebra
$\Sigm^{\roman c}[M]$
on $M$
has the form
$$
\Sigm^{\roman c}[M] = \oplus \Sigm_k^{\roman c}[M] 
$$
where $\Sigm_0^{\roman c}[M]=R$ and where, for $k \geq 1$,
$\Sigm_k^{\roman c}[M]$ is the differential graded subspace
$$
\Sigm_k^{\roman c}[M] =(M^{\otimes k})^{\Sigma_k}
\subseteq \roman T^{\roman c}_k[M] = M^{\otimes k}
$$
of invariants under the canonical action 
(in the graded sense) of
the symmetric group $\Sigma_k$ on $k$-letters
on the $k$'th tensor power
$M^{\otimes k}$;
see e.~g. \cite\extensta.

Let  $\fra g$
be a chain complex, endowed with a graded skew-symmetric 
bracket $[\ ,\ ] \colon \fra g \otimes \fra g \to \fra g$
of degree zero, not necessarily a graded Lie bracket
nor necessarily a chain map.
Consider the differential graded symmetric coalgebra
$\Sigm^{\roman c}[s \fra g]$,
the differential $d$ 
on $\Sigm^{\roman c}[s \fra g]$
being induced from that 
on $\fra g$,
and let $\partial$ be the coderivation
$\Sigm^{\roman c}[s \fra g] \to \Sigm^{\roman c}[s \fra g]$
of degree $-1$ induced by 
$[\ ,\ ]$;
it is manifest that
$\partial$
lowers the  obvious filtration of
$\Sigm^{\roman c}[s \fra g]$.
We then have the following folk-lore result a proof of which is straightforward
and left to the reader.

\proclaim{Proposition 1.3}
The bracket $[\ ,\ ]$
turns 
$\fra g$ into
a dg Lie algebra 
if and only if
$\partial$ is a coalgebra
perturbation 
of the differential $d$.
Furthermore,
any differential graded Lie bracket
on $\fra g$
is determined by the coalgebra
perturbation 
induced from the bracket. \qed
\endproclaim

The more \lq sophisticated\rq\ description
of a dg Lie algebra alluded to earlier is that
in terms of the
coalgebra perturbation $\partial$ given in Proposition 1.3.
Below we will abstract from this description and take
coalgebra perturbations which do not necessarily
correspond to true Lie brackets.
Suffice it to explain that,
under the circumstances of Proposition 1.3, the vanishing
of $(d + \partial)^2$ is equivalent 
to the bracket being a chain map satisfying 
the graded Jacobi identity.
More precisely,
the vanishing
of $(d + \partial)^2$ is equivalent 
to that of
$d\partial + \partial d$
and
$\partial \partial$ separately;
now
$d\partial + \partial d$
to be zero corresponds to the bracket being a chain map,
and
the vanishing of 
$\partial \partial$
is equivalent to the graded Jacobi identity.
When  $\fra g$
is an ordinary Lie algebra
over a field $\fiel k$,
$\Sigm^{\roman c}[s \fra g]$
with this differential $\partial$
boils down to the ordinary Koszul or
Chevalley-Eilenberg complex calculating the homology of $\fra g$
with coefficients in $\fiel k$;
the dual
$\roman{Hom}(\Sigm^{\roman c}_{[\cdot,\cdot]}[s\fra g],\fiel k)$
is then the Chevalley-Eilenberg algebra
calculating the cohomology of $\fra g$ with coefficients in $\fiel k$.
For a general differential graded Lie algebra $\fra g$,
with Lie bracket $[\ ,\ ]$,
we will occasionally
write
$\Sigm^{\roman c}_{[\cdot,\cdot]}[s\fra g]$
for the resulting
coaugmented  differential graded cocommutative coalgebra
or,
more simply,
following \cite\quilltwo\ (Appendix) just $\Cal C [\fra g]$
and refer to it
as its {\it generalized Koszul\/} or {\it Chevalley-Eilenberg complex\/}.

We now explain the 
{\it Master Equation\/}
and the formal circumstances thereof.
In general, the master equation  makes sense
as an equation on elements of a differential graded associative algebra
$\Cal A$;
it then has the form
$$
D \Gammat = \Gammat \Gammat
\tag1.4.1
$$
where
$\Gammat$ is a homogeneous element of $\Cal A$
(necessarily of degree $-1$), where $D$ refers to the differential
in $\Cal A$, and where
$\Gammat \Gammat$ is the product in $\Cal A$ of
$\Gammat$ with itself.
Likewise,
as an equation on elements of a differential graded Lie algebra
$\Cal L$,
for an element $\Gammat$  of $\Cal L$,
the master equation has the form
$$
D \Gammat = \frac 12 [\Gammat, \Gammat]. 
\tag1.4.2
$$
When $\Cal A$ is the universal differential graded algebra
of a differential graded Lie algebra
$\Cal L$,
for an element $\Gammat$  of $\Cal L$,
the master equations (1.4.1) and (1.4.2)
are manifestly equivalent.
Inspection shows that,
given
a differential graded algebra
$\Cal A$ and a solution $\Gammat$ of the master equation,
the operator $d_{\Gammat}$ on $\Cal A$
defined by
$$
d_{\Gammat}
(a) = d a - \Gammat a,\quad
a \in \Cal A,
\tag1.4.3
$$
yields a new
differential graded algebra structure on $\Cal A$.
Likewise,
given
a differential graded Lie algebra
$\Cal L$ and a solution $\Gammat$ of the master equation,
the operator $d_{\Gammat}$ on $\Cal L$
defined by
$$
d_{\Gammat}
(a) = d a - [\Gammat, a],\quad
a \in \Cal L,
\tag1.4.4
$$
yields a new
differential graded Lie algebra structure on $\Cal L$.

\smallskip\noindent
{\smc Remark.}
In the literature,
instead of (1.4.2),
it is also customary to take
the equation
$D \Gammat + \frac 12 [\Gammat, \Gammat] = 0$
as (Lie algebra) master equation. 
The 
formula (1.4.4) for the new operator 
must then accordingly be replaced by
$d_{\Gammat} (a) = d a + [\Gammat, a]$.

Let  $\fra g$
be a 
differential graded Lie algebra.
We then have the differential $d$ on
$\Cal C [\fra g]=\Sigm^{\roman c}[s \fra g]$
induced from that on $\fra g$
and, cf. Proposition 1.3,
the coalgebra perturbation
thereof induced from the graded Lie structure on $\fra g$.
More general
coalgebra perturbations
of the differential $d$ on
$\Sigm^{\roman c}[s \fra g]$
are related to
solutions of a certain master equation 
as we will explain in Section 2. 
This master equation involves
a differential graded Lie algebra 
(different from $\fra g$)
of a kind which we now explain.
Recall that, for any coaugmented graded coalgebra $C$
and any graded Lie algebra $\fra h$,
given $a,b \colon C \to \fra h$, 
(with a slight abuse of notation)
their {\it cup bracket\/}
$[a, b]$ is given by the composite
$$
C @>{\Delta}>> C\otimes C @>{a\otimes b}>>\fra h \otimes\fra h @>
{[\cdot,\cdot]}>> \fra h.
\tag1.5.1
$$
Now
coalgebra perturbations
of the differential $d$ on
$\Cal C [\fra g]=\Sigm^{\roman c}[s \fra g]$
(which is induced from that on $\fra g$)
are related to solutions of the master equation
in
$\Cal L = \roman{Hom}(C, \fra g)$
where $C$ is an appropriate
 differential graded cocommutative coalgebra
sufficiently closely related to $\fra g$
(e.~g. $C=\Sigm^{\roman c}[s \fra g]$, 
endowed with the differential 
induced from that on $\fra g$).
Such solutions of the master equation 
may be described in the language of twisting cochains.
This is crucial for 
our method 
aimed at constructing solutions of the master equation
which we will explain in Section 2 below.

Recall that,
given
a general differential graded
coalgebra $C$ and
a general differential graded
algebra $A$,
for $a,b \colon C \to A$, their {\it cup product\/}
$a \smile b$ is given by the composite
$$
C @>{\Delta}>> C \otimes C @>{a \otimes b}>> A \otimes A @>\mu>> A
\tag1.5.2
$$
where
$\Delta$ and $\mu$ denote the 
structure maps.
This cup product turns $\roman{Hom}(C,A)$
into a differential graded algebra.                               
Furthermore, a coaugmentation  and augmentation
map
of $C$ and $A$ respectively induce an augmentation map
for
$\roman{Hom}(C,A)$.
\smallskip
\noindent
{\smc Definition 1.5.}
Given
a coaugmented differential graded
coalgebra $C$ and
an augmented differential graded
algebra $A$,
a {\it twisting cochain\/}
$t \colon C \to A$ is a homogeneous morphism of
degree $-1$
whose composites with the coaugmentation and augmentation maps are zero
and which
satisfies
$$
Dt = t \smile t.
\tag1.5.3
$$
Given
a  coaugmented differential graded cocommutative
coalgebra $C$ and
a differential graded Lie algebra $\fra h$,
a {\it Lie algebra twisting cochain\/}
$t \colon C \to \fra h$ is a homogeneous morphism of
degree $-1$
whose composite with the coaugmentation  map is zero
and which
satisfies
$$
Dt = \frac 12 [t,t];
\tag1.5.4
$$ 
cf. \cite\moorefiv,
\cite\quilltwo.

In particular, when $A$ is the universal differential
graded algebra of a 
differential graded
Lie algebra $\fra h$
and when $C$ is graded cocommutative,
given a morphism 
$t \colon C \to A$
of degree $-1$
whose values 
lie in $\fra h$, the defining property 
(1.5.3)
for $t$ being a twisting cochain
is equivalent 
to the defining property 
(1.5.4)
for
$t$ being a Lie algebra twisting cochain.

For later reference we recall that,
for any differential graded Lie algebra $\fra h$,
its universal Lie algebra
twisting cochain
$$
\tau_{\fra h} \colon 
\Cal C [\fra h]= \Sigm^{\roman c}_{[\ ,\ ]}[s\fra h]
@>>>
\fra h,
\tag1.5.5
$$
cf. e.~g. Appendix B of \cite\quilltwo,
is given by $\tau_{\fra h}(sx) = x$ for $x \in \fra h$
and
$\tau_{\fra h}(y) = 0$ when $y \in \Sigm^{\roman c}_k[s\fra h]$ for $ k \ne 1$.

\medskip\noindent{\bf 2. Homological Perturbation Theory (HPT)}
\smallskip\noindent
Homological perturbation theory is concerned with transferring
various kinds of algebraic structure through a homotopy equivalence.
We first describe an essential piece of machinery.
\smallskip\noindent
{\smc Definition 2.1.}
A {\it contraction\/}
$$
\Nsddata N{\nabla}{\pi}Mh
\tag 2.1.1
$$
of chain complexes,
referred to as well as
SDR-{\it data\/} in the literature,
consists of
\smallskip
-- chain complexes $N$ and $M$,
\newline
\indent
-- chain maps 
$\pi\colon N \to M$ and $\nabla \colon M \to N$,
\newline
\indent
--  a morphism $h\colon N \to N$ 
of the underlying graded modules of degree 1;
\newline
\noindent
these data are required to satisfy
$$
\align
\pi \nabla &= \roman{Id}
\tag2.1.2
\\
Dh &= \nabla \pi - \roman{Id}
\tag2.1.3
\\
\pi h &= 0, \quad h \nabla = 0,\quad hh = 0. 
\tag2.1.4
\endalign
$$
The requirements (2.1.4) are referred to as
{\it annihilation properties\/}
or {\it side conditions\/};
they can always be achieved without loss of generality,
cf. e.~g. what is said under \lq Remarks\rq\ 
in Section 1 of \cite\gugenhth.

Given a contraction (2.1.1),
we will say that
$N$ {\it contracts onto\/} $M$.
If, furthermore, $N$ and $M$ are filtered chain complexes, and if 
$\pi$, $\nabla$ and $h$ are filtration preserving,
the contraction is said to be {\it filtered\/}.
This notion of contraction was introduced in \S 12 of
\cite\eilmactw;
it is among the basic notions in homological perturbation theory,
cf. \cite{\gulstatw,\, \huebkade}
and the literature there.

Let at first $\fra g$ be a chain complex, and let
$$
\Nsddata {\fra g} {\nabla}{\pi}{\roman H(\fra g)}h
\tag2.2
$$
be a contraction of chain 
complexes.
Its existence would, in general, be an additional assumption
but such a contraction
will always exist when
$\roman H(\fra g)$ is free over the ground ring $R$ or when 
$R$ is a field.
Under our circumstances,
since we assume that the ground ring $R$ contain the rationals
as a subring,
the existence of a contraction is a very mild assumption, if any.
Notice that we momentarily ignore additional structure on $\fra g$.

Consider the induced {\it filtered\/} contraction
$$
\Nsddata {\Sigm^{\roman c}[s\fra g]} {\Sigm^{\roman c}\nabla}
{\Sigm^{\roman c}\pi}
{\Sigm^{\roman c}[s\roman H( \fra g)]}{\Sigm^{\roman c}h}
\tag2.3
$$
of {\it differential graded coaugmented coalgebras\/},
the filtrations being the ordinary coaugmentation filtrations.
It may be obtained in the following way:

Any contraction
$\Nsddata N{\nabla}{\pi}Mh$
of chain complexes
induces a filtered contraction
$$
\Nsddata 
{\roman T^{\roman c}[N]}
{\roman T^{\roman c}\nabla}
{\roman T^{\roman c}\pi}
{\roman T^{\roman c}[M]}
{\roman T^{\roman c}h}
\tag2.4
$$
of coaugmented differential graded coalgebras.
A version thereof 
is spelled out as a {\it contraction of bar constructions\/}
already in Theorem 12.1 of 
\cite\eilmactw; 
the filtered contraction (2.4) may be found in
\cite\gugenlam\ (3.2), \cite\gulasta\ (\S 3),
\cite\gulstatw\  (2.2);
the dual filtered contraction
of augmented differential graded algebras
which involves the differential graded
tensor algebra
has been given in
\cite\huebkade\ ${\roman (2.2.0)}^*$. 
The differential graded
symmetric coalgebras
${\Sigm^{\roman c}[M]}$
and
${\Sigm^{\roman c}[N]}$
being differential graded subcoalgebras of
${\roman T^{\roman c}[M]}$
and
${\roman T^{\roman c}[N]}$, respectively,
the morphisms
${\roman T^{\roman c}\nabla}$ and
${\roman T^{\roman c}\pi}$ 
pass to corresponding morphisms
${\Sigm^{\roman c}\nabla}$
and
${\Sigm^{\roman c}\pi}$ 
respectively, 
and
${\Sigm^{\roman c}h}$
arises from
${\roman T^{\roman c}h}$
by {\it symmetrization\/}, so that
$$
\Nsddata 
{\Sigm^{\roman c}[N]}
{\Sigm^{\roman c}\nabla}
{\Sigm^{\roman c}\pi}
{\Sigm^{\roman c}[M]}
{\Sigm^{\roman c}h}
\tag2.5
$$
constitutes a filtered
contraction of coaugmented differential graded coalgebras.
Here 
$\Sigm^{\roman c}\nabla$
and
$\Sigm^{\roman c}\pi$ are morphisms of differential graded coalgebras
but, beware,
even though
${\roman T^{\roman c}h}$
is compatible with the coalgebra structure in the sense that
it is a homotopy of morphisms of differential graded coalgebras,
$\Sigm^{\roman c}h$ no longer has such a compatibility property
in a naive fashion.
Indeed, for cocommutative coalgebras, the notion of homotopy
is a subtle concept, cf. \cite\schlstas.

\smallskip\noindent
{\smc Definition 2.6.}
Given a chain complex $\fra g$,
an sh-{\it Lie structure\/} or $L_{\infty}$-{\it structure\/}
on $\fra g$ is
a {\it coalgebra perturbation\/} $\partial$
of the differential $d$ on
the coaugmented differential graded symmetric coalgebra
$\Sigm^{\roman c}[s\fra g]$ on $s\fra g$, that is,
an operator $\partial$ of degree $-1$
which is compatible with the 
coalgebra structure,
lowers filtration by 1, and
satisfies
$$
(d + \partial)^2 = 0
\tag2.6.1
$$
or, equivalently,
$$
[d,\partial] + \partial \partial = 0
\tag2.6.2
$$
and, furthermore,
$$
\partial \eta = 0
\tag2.6.3
$$
so that
the sum $d + \partial$ endows
$\Sigm^{\roman c}[s\fra g]$
with a new coaugmented differential graded  coalgebra structure.

In contrast to the case of a strict dg Lie algebra, here we may
have non-zero
terms
$\Sigm^{\roman c}_k[s\fra g] \to s\fra g$
for $k \geq 2$. For example, the Jacobi identity may hold only
modulo an exact term provided by
a morphism of the kind
$\Sigm^{\roman c}_3[s\fra g] \to s\fra g$.

Given an sh-Lie structure $\partial$ on $\fra g$,
we write
$\Sigm^{\roman c}_{\partial}[s\fra g]$
for the new coaugmented differential graded coalgebra.
We have explained in Proposition 1.3
how a graded Lie bracket 
$[\cdot,\cdot]$
on $\fra g$
yields such a coalgebra perturbation
$\partial$;
the resulting chain complex
$\Sigm^{\roman c}_{\partial}[s\fra g]$
then boils down to
the generalized Koszul complex  denoted there by 
$\Sigm^{\roman c}_{[\cdot,\cdot]}[s\fra g]$
or just
$\Cal C [\fra g]$.

Given two
sh-Lie algebras 
$(\fra g_1, \partial_1)$ and
$(\fra g_2, \partial_2)$,
an {\it sh-morphism\/}
or
{\it sh-Lie map\/}
from
$(\fra g_1, \partial_1)$ to
$(\fra g_2, \partial_2)$
is
a morphism
${
\Sigm^{\roman c}_{\partial_1}[s\fra g_1]
@>>>
\Sigm^{\roman c}_{\partial_2}[s\fra g_2]
}$
of
differential graded coalgebras.

\proclaim{Theorem 2.7}
Given a differential graded Lie algebra
$\fra g$
and a contraction
of chain complexes of the kind {\rm (2.2)},
the corresponding coalgebra perturbation  of the differential
on $\Sigm^{\roman c}[s\fra g]$ being written $\partial$,
the data determine
\newline\noindent
{\rm (i)}
a  differential $\Cal D$
on 
$\Sigm^{\roman c}[s\roman H( \fra g)]$
turning the latter into a coaugmented differential graded coalgebra 
(i.~e. $\Cal D$ is a coalgebra perturbation  of the zero differential)
and hence endowing
$\roman H( \fra g)$
with an sh-Lie algebra structure
and 
\newline\noindent
{\rm (ii)}
a Lie algebra twisting cochain $\tau \colon 
\Sigm_{\Cal D}^{\roman c}[s\roman H( \fra g)] \to \fra g$
whose adjoint $\overline \tau$, written
\linebreak
$(\Sigm^{\roman c}\nabla)_{\partial} \colon 
\Sigm_{\Cal D}^{\roman c}[s\roman H( \fra g)] \to \Cal C[\fra g]$,
induces an isomorphism on homology.
\newline\noindent
Furthermore, 
$(\Sigm^{\roman c}\nabla)_{\partial}$
admits an extension to a new contraction
$$
\Nsddata {\Sigm^{\roman c}_{\partial}[s\fra g]} 
{(\Sigm^{\roman c}\nabla)_{\partial}}
{(\Sigm^{\roman c}\pi)_{\partial}}
{\Sigm^{\roman c}_{\Cal D}[s\roman H( \fra g)]}{(\Sigm^{\roman c}h)_{\partial}}
\tag2.7.1
$$
of filtered chain complexes (not nececessarily of coalgebras).
\endproclaim

For intelligibility we point out that,
according to the convention introduced
before,
in the statement of the theorem,
the perturbed
differential graded coalgebra 
which corresponds to the 
asserted sh-Lie algebra structure
on $\roman H( \fra g)$
is written
$\Sigm^{\roman c}_{\Cal D}[s\roman H( \fra g)]$.
Further,
the morphisms 
$(\Sigm^{\roman c}\nabla)_{\partial}$
and
$(\Sigm^{\roman c}\pi)_{\partial}$
coming into play 
in the new contraction (2.7.1)
are sh-Lie maps between
$(\roman H(\fra g),\Cal D)$ and $(\fra g,\partial)$.
They are even
sh-equivalences (quasi-isomorphisms) in the sense that they induce
isomorphisms on homology.
It is worthwhile pointing out that the induced bracket on
$\roman H(\fra g)$
is a strict graded Lie bracket but, in general,
the differential $\Cal D$ involves
meaningful terms of higher order.

\demo{Sketch of the proof of Theorem 2.7} For $b \geq 1$,
write $\Sigm_b^{\roman c}$
for the homogeneous degree $b$ component of
$\Sigm^{\roman c}[s\roman H( \fra g)]$.
The operator $\Cal D$ and twisting cochain $\tau$ are obtained
as infinite series by the following recursive procedure where $b \geq 2$:
$$
\alignat 1
\tau &= \tau^1 + \tau^2 + \dots, \quad 
\quad \tau^1 = \nabla \tau_{\roman H (\fra g)},
\quad \tau^j \colon\Sigm_j^{\roman c} \to \fra g,\ j \geq 1,
\tag2.7.2
\\
\tau^b &= \frac 12 h([\tau^1,\tau^{b-1}] +  \dots + [\tau^{b-1},\tau^1])
\\
\Cal D &= \Cal D^1 + \Cal D^2 + \dots
\tag2.7.3
\endalignat
$$
where $\Cal D^{b-1}$
is
the coderivation 
of $\Sigm^{\roman c}[s\roman H( \fra g)]$
determined by
$$
\tau_{\roman H(\fra g)} \Cal D^{b-1} =
\frac 12 \pi
([\tau^1,\tau^{b-1}] +  \dots + [\tau^{b-1},\tau^1])
\colon
\Sigm_b^{\roman c}
\to \roman H (\fra g).
$$
The sums (2.7.2) and (2.7.3) are infinite but, applied to a specific element
which, in view of the assumptions,
necessarily lies in some
finite filtration degree subspace,
only finitely many terms will be non-zero,
whence the convergence is naive.
For example, cf. (2.7.3),
for $n \geq 1$,
the operator $\Cal D^n$ vanishes on
the constituent
$F_n$ 
of the coaugmentation filtration 
$$
R=F_0\subseteq F_1 \subseteq F_2 \subseteq \dots
$$
of
$\Sigm^{\roman c}[s\roman H( \fra g)]$.
Since the coaugmentation filtration is cocomplete, 
if the infinite sum 
(2.7.3)
for $\Cal D$ is applied to a particular
element,
only finitely many terms are non-zero. 
The summand $\Cal D^1$ is the ordinary Cartan-Chevalley-Eilenberg
differential for the classifying coalgebra $\Cal C[\roman H(\fra g)]$
of the graded Lie algebra $\roman H(\fra g)$.
It may be shown that
$\Cal D$ is indeed a coalgebra differential and that $\tau$ is a twisting 
cochain; the details will be given elsewhere.
A construction which is formally similar but yet substantially different
may be found in Proposition 2.1 of \cite\gugenhth.

A spectral sequence
comparison argument shows that the adjoint
$\overline \tau =(\Sigm^{\roman c}\nabla)_{\partial}$
of $\tau$
induces an isomorphism on homology.
Hence $(\Sigm^{\roman c}\nabla)_{\partial}$
admits an extension to a contraction
of chain complexes 
of the kind (2.7.1). \qed
\enddemo

\smallskip
\noindent{\smc 2.8. Refinements.}

\proclaim{Addendum 2.8.1}
Under the circumstances of Theorem {\rm 2.7},
$$
\tau \colon
{\Sigm^{\roman c}_{\Cal D}[s\roman H( \fra g)]}
@>>>
\fra g,
\tag2.8.2
$$
viewed as an element of degree $-1$ 
of 
the differential graded Lie algebra
\linebreak
$\roman{Hom}(\Sigm^{\roman c}_{\Cal D}[s\roman H( \fra g)], \fra g)$,
satisfies the master equation
{\rm (1.4.2)}.
\endproclaim

Here
$
\roman{Hom}(\Sigm^{\roman c}_{\Cal D}[s\roman H( \fra g)], \fra g)
$
is endowed with the graded cup bracket (1.5.1)
induced by the graded coalgebra structure
on
$\Sigm^{\roman c}_{\Cal D}[s\roman H( \fra g)]$ and the graded
bracket on $\fra g$.
The twisting cochain (2.8.2)
is our {\it most general solution of the master
equation\/};
the other (more special) solutions of the master equation derive from it.

\proclaim{Addendum 2.8.3}
Under the circumstances of Theorem {\rm 2.7},
suppose in addition that there is
a  differential $\widetilde {\Cal D}$
on 
$\Sigm^{\roman c}[s\roman H( \fra g)]$
turning the latter into a coaugmented differential graded coalgebra
in such a way that
${
(\Sigm^{\roman c} \pi) \partial = \widetilde {\Cal D}(\Sigm^{\roman c} \pi) .
}$
Then
$\Cal D =\widetilde {\Cal D}$
and
$(\Sigm^{\roman c}\pi)_{\partial}$
may be taken to be $\Sigm^{\roman c}\pi$.
In particular,
when
$(\Sigm^{\roman c} \pi) \partial$ is zero,
the differential $\Cal D$ on
$\Sigm^{\roman c}[s\roman H( \fra g)]$
is necessarily zero, that is,
the new contraction {\rm (2.7.1)} has the form
$$
\Nsddata {\Sigm^{\roman c}_{\partial}[s\fra g]} 
{(\Sigm^{\roman c}\nabla)_{\partial}}
{\phantom{xy} \Sigm^{\roman c}\pi \phantom{xy}}
{\Sigm^{\roman c}[s\roman H( \fra g)]}{(\Sigm^{\roman c}h)_{\partial}}.
\tag2.8.4
$$
For example, this happens to be the case
when the composite 
${ \fra g \otimes \fra g
@>{[\cdot,\cdot]}>>
\fra g
@>\pi>> \roman H(\fra g)
}$
is zero.
\endproclaim

In fact, this follows at once from the descriptions
(2.7.3) and (2.7.2)
for   $\Cal D$ and $\tau$:
Since
${
(\Sigm^{\roman c} \pi) \partial = \widetilde {\Cal D}(\Sigm^{\roman c} \pi)
}$,
only the quadratic part 
$\widetilde {\Cal D}^1 \colon 
\Sigm_2^{\roman c} \to s\roman H(\fra g)$
of
$\widetilde {\Cal D}$
is non-zero and this quadratic part
necessarily coincides
with the quadratic part 
$\Cal D^1$
of
$\Cal D$.
Furthermore,
since the constituent $\pi$ of the contraction (2.2) is necessarily
a morphism of differential graded Lie algebras,
$\Cal D$ has no higher terms as well.
Indeed, since $\pi h=0$, for $b \geq 2$, $\pi\tau^b = 0$
whence 
$\pi [\tau^j,\tau^k] =[\pi \tau^j,\pi\tau^k]$
($j,k \geq 1$)
is non-zero only for $j=k=1$.
Hence 
$\Cal D^b$ is zero for $b \geq 2$.

Likewise, we have the following.

\proclaim{Addendum 2.8.5}
Under the circumstances of Theorem {\rm 2.7},
suppose in addition that there is
a  differential $\widetilde {\Cal D}$
on 
$\Sigm^{\roman c}[s\roman H( \fra g)]$
turning the latter into a coaugmented differential graded coalgebra
in such a way that
${
\partial (\Sigm^{\roman c} \nabla)  = 
(\Sigm^{\roman c} \nabla)\widetilde {\Cal D} .
}$
Then
$\Cal D =\widetilde {\Cal D}$
and
$(\Sigm^{\roman c}\nabla)_{\partial}=\Sigm^{\roman c}\nabla$.
In particular,
when
$ \partial (\Sigm^{\roman c} \nabla)$ is zero,
the differential $\Cal D$ on
$\Sigm^{\roman c}[s\roman H( \fra g)]$
is necessarily zero, that is,
the new contraction {\rm (2.7.1)} has the form
$$
\Nsddata {\Sigm^{\roman c}_{\partial}[s\fra g]} 
{\phantom{xy}\Sigm^{\roman c}\nabla \phantom{}}
{(\Sigm^{\roman c}\pi)_{\partial}}
{\Sigm^{\roman c}[s\roman H( \fra g)]}{(\Sigm^{\roman c}h)_{\partial}}.
\tag2.8.6
$$
For example, this happens to be the case
when the composite 
\linebreak
${ 
\roman H\fra g \otimes \roman H\fra g
@>{\nabla \otimes \nabla}>>
\fra g \otimes \fra g
@>{[\cdot,\cdot]}>> \fra g
}$
is zero.
\endproclaim

In the case of (2.8.3),
$\pi\colon \fra g \to \roman H(\fra g)$ is in fact a strict
graded Lie morphism 
and so is, likewise, 
$\nabla \colon \roman H(\fra g)\to \fra g$ in the case of (2.8.5).
We note that observations similar to these addenda
 may be found in
(2.4.1) and (2.4.2) of \cite\perturba.

Under the circumstances of
Theorem 2.7 where 
$\fra g$ is a true differential graded
Lie algebra,
the perturbation $\Cal D$ still endows $\roman H( \fra g)$
with an  sh-Lie algebra structure
which is in general non-trivial,
and this will be so even if the induced Lie bracket
on
$\roman H( \fra g)$ is zero
but {\it not\/} if an additional 
(technical) 
condition 
(which arises from abstraction from the statement of the 
$\partial \overline \partial$-Lemma for K\"ahler manifolds,
see what is said below)
is satisfied,
according to Addendum 2.8.3.
This latter is the situation for Barannikov-Kontsevich \cite\barakont\ 
and Manin \cite\maninbtw\ (III.10),\,\cite\maninfiv\ (Section 6).
They have stronger conditions in the main applications,
for which we spell out the following.

\proclaim{Theorem 2.9}
Given a differential graded Lie algebra $\fra g$,
a differential graded 
Lie subalgebra $\fra m$ of $\fra g$,
and a contraction
$$
\Nsddata {\fra m} {\nabla}{\pi}{\roman H(\fra g)}h
\tag2.9.1
$$
of chain complexes so that
the composite
${ \fra m \otimes \fra m
@>{[\cdot,\cdot]}>>
\fra m
@>\pi>> \roman H(\fra g)
}$
is zero,
the induced bracket on $\roman H( \fra g)$ is zero,
that is, as a graded Lie algebra, $\roman H( \fra g)$ is abelian,
and
the  data determine a solution 
${
\tau \in \roman{Hom}(\Sigm^{\roman c}[s\roman H( \fra g)], \fra g)
}$
of the master equation
$$
d \tau  =\frac 12 [\tau,\tau]
\tag2.9.2
$$
in such a way that 
the following hold:
\newline\noindent
{\rm (i)} 
The composite $\pi \tau$ coincides with the universal twisting cochain
$\Sigm^{\roman c}[s\roman H( \fra g)] \to  \roman H( \fra g)$
for the abelian graded Lie algebra $\roman H( \fra g)$;
\newline\noindent
{\rm (ii)} 
the values of $\tau$ lie in $\fra m$. 
\endproclaim

Notice that (2.9.2) is somewhat simpler than 
the general master equation (1.4.2)
since there is no non-zero differential on
$\Sigm^{\roman c}[s\roman H( \fra g)]$.

\demo{Proof}
Under these circumstances,
the homology of $\fra m$ is necessarily 
identified with $\roman H(\fra g)$
under (2.9.1).
In view of Addendum 2.8.3,
the corresponding contraction 
(2.7.1)
has the form
$$
\Nsddata {\Sigm^{\roman c}_{\partial}[s\fra m]} 
{(\Sigm^{\roman c}\nabla)_{\partial}}
{\phantom{xy}\Sigm^{\roman c}\pi\phantom{z}}
{\Sigm^{\roman c}[s\roman H( \fra g)]}{(\Sigm^{\roman c}h)_{\partial}}.
$$
Notice that
the differential on $\Sigm^{\roman c}[s\roman H( \fra g)]$
is zero.
Consider
the corresponding
Lie algebra twisting cochain (2.8.2) 
(where the role of $\fra g$ in (2.8.2) is now played by
$\fra m$) which we write as
$$
\widetilde \tau =
\tau_{\fra m} \circ (\Sigm^{\roman c}\nabla)_{\partial}
\colon
{\Sigm^{\roman c}[s\roman H( \fra g)]}
@>>>
\fra m,
$$
and denote by $\tau$ its composite with the inclusion
$\fra m \subseteq\fra g$.
Since this inclusion is a morphism of differential graded Lie
algebras,
$\tau$ has the asserted properties.  \qed
\enddemo

%\beginsection 3. Differential Gerstenhaber and Batalin-Vilkovisky algebras
\medskip\noindent
{\bf 3. Differential Gerstenhaber and Batalin-Vilkovisky algebras}
\smallskip\noindent
Recall that a {\it Gerstenhaber algebra\/}
is a graded commutative $R$-algebra
$\Cal A$
together with a graded Lie bracket
from $\Cal A \otimes_R \Cal A$ to $\Cal A$
of degree $-1$ (in the sense that,
if $\Cal A$ is regraded down by one,
$[\cdot,\cdot]$ is an ordinary graded Lie bracket)
such that, for each homogeneous element $a$ of $\Cal A$,
the operation
$[a,\cdot]$ is a derivation of $\Cal A$ of degree $|a|-1$
where $|a|$ refers to the degree of $a$;
see \cite\geschthr\ 
where these objects are called G-algebras,
or \cite{\kosmathr,\,\liazutwo,\,\xuone};
for a Gerstenhaber algebra
$\Cal A$,
the bracket
from $\Cal A \otimes_R \Cal A$ to $\Cal A$
will henceforth be referred to as its {\it Gerstenhaber\/}
bracket.

\smallskip\noindent
{\smc Definition 3.1.}
A {\it differential Gerstenhaber algebra\/}
$(\Cal A, [\cdot,\cdot], d)$
consists of a
Gerstenhaber algebra $(\Cal A, [\cdot,\cdot])$
and a differential $d$ of degree $+1$ on $\Cal A$ which is a derivation
for 
the multiplication of $\Cal A$;
$(\Cal A, [\cdot,\cdot], d)$ will be said to be a {\it strict\/} 
differential Gerstenhaber algebra
provided  
the differential $d$ is a derivation
for the
Gerstenhaber bracket
$[\cdot,\cdot]$ as well, that is, 
$$
d[x,y] = [dx,y] -(-1)^{|x|} [x,dy], \quad x,y \in \Cal A.
\tag3.1.1
$$

Let
$(\Cal A, [\cdot,\cdot], d)$
be a
strict differential
Gerstenhaber algebra.
We denote the corresponding ordinary differential graded Lie algebra,
with differential spelled out as a homology differential, by
$\fra g^{\Cal A}$ or, more simply,
by $\fra g$, when there is no risk of confusion.
Thus,
$\fra g_{-*} = \Cal A^{*+1}$, that is to say,
$$
\fra g_1 = \Cal A^0,
\quad
\fra g_0 = \Cal A^1,
\quad
\fra g_{-1} = \Cal A^2,
\quad
\dots ,
\quad
\fra g_{-n} = \Cal A^{n+1},
\tag3.2
$$
so that
the graded bracket $[\cdot,\cdot]$
and differential $d$
on $\fra g$ are
of the ordinary kind,
i.~e.  of the form
$$
[\cdot,\cdot]\colon
\fra g_j \otimes \fra g_k @>>>
\fra g_{j+k},
\quad
d \colon \fra g_j @>>> \fra g_{j-1}.
$$
Notice that,
when we write
$\fra a_{-*} = \Cal A^*$, so that
$\fra a_{0} = \Cal A^0$,
$\fra a_{-1} = \Cal A^1$, etc., we have
$\fra g = s \fra a$ and
$\fra a = s^{-1} \fra g$.
The notation $\fra a$ is just used to rewrite
the strict differential Gerstenhaber algebra
$\Cal A$ in homology degrees.

Consider a contraction
of the kind (2.2).
Application of Theorem 2.7,
the operator $\partial$ 
on
$\Sigm^{\roman c}_{\partial}[s\fra g]$
which corresponds to the
Lie bracket on $\fra g$
being regarded as the perturbation to be transferred
(to $\Sigm^{\roman c}[s\roman H( \fra g)]$),
yields a new contraction 
$$
\Nsddata {\Sigm^{\roman c}_{\partial}[s\fra g]} 
{(\Sigm^{\roman c}\nabla)_{\partial}}
{(\Sigm^{\roman c}\pi)_{\partial}}
{\Sigm^{\roman c}_{\Cal D}[s\roman H( \fra g)]}{(\Sigm^{\roman c}h)_{\partial}}
$$
of filtered differential graded coalgebras.
The 
corresponding 
Lie algebra twisting cochain denoted 
in  (2.8.2) 
by $\tau$,
now written out as having $\fra a$
as its target rather than $\fra g$,
appears as an element
${
\tau \in
\roman{Hom}(\Sigm^{\roman c}_{\Cal D}[s\roman H( \fra g)], \fra a)
}$
of degree $-2$ 
satisfying
the master equation
$$
D \tau = \frac 12 [\tau,\tau].
\tag 3.3
$$
Here $D$ is the Hom-differential,
and $[\cdot,\cdot]$
refers to the graded cup bracket 
(1.5.1)
which is induced
by the graded coalgebra structure
on $\Sigm^{\roman c}[s\roman H( \fra g)]$
and the graded Lie structure on $\fra g$ which, in turn,
is the graded Gerstenhaber bracket on $\fra a$.

For consistency with what is in the literature,
we now rewrite 
$\roman{Hom}(\Sigm^{\roman c}_{\Cal D}[s\roman H( \fra g)], \fra a)$
and the element
$\tau$ thereof
in terms of cohomology degrees:
Let
$\Cal G^* = \fra g_*$;
in this description,
$\Cal A = S\Cal G$, 
where
$S$ refers to the suspension operator in {\it cohomology degrees\/}, so that,
for every integer $j$,
$$
(\Cal A)^j=(S\Cal G)^j = \Cal G^{j+1},
\quad
(S^{-1}\roman H^*(\Cal G))^j
=
\roman H^{j+1}(\Cal G).
$$
The twisting cochain
$\tau$
now appears as an element
$$
\tau \in
\Cal B =\roman{Hom}(\Sigm^{\roman c}_{\Cal D}[S^{-1}\roman H(\Cal G)], \Cal A)
\tag3.4
$$
of degree $+2$ 
satisfying the master equation.

\smallskip
\noindent
{\smc Definition 3.5.}
For a Gerstenhaber algebra $\Cal A$
over $R$,
with bracket operation written $[\cdot,\cdot]$,
an $R$-linear operator $\Delta$
on $\Cal A$  
of degree $-1$
is said to
{\it generate\/}
the Gerstenhaber bracket
provided, for every homogeneous $a, b \in \Cal A$,
$$
[a,b] = (-1)^{|a|}\left(
\Delta(ab) -(\Delta a) b - (-1)^{|a|} a (\Delta b)\right);
\tag3.5.1
$$
the operator $\Delta$ is then called a {\it generator\/}.
A generator $\Delta$  is said to be {\it exact\/}
provided $\Delta\Delta$ is zero, that is, $\Delta$ is a differential.
A Gerstenhaber algebra $\Cal A$
together with a generator $\Delta$
will be called a {\it weak Batalin-Vilkovisky\/} algebra
(or weak BV-algebra);
when
the generator is exact,
$(\Cal A,\partia)$
is (more simply)
called a
{\it Batalin-Vilkovisky\/} algebra
(or BV-algebra).
The notation $\Delta$ for a generator of a Batalin-Vilkovisky
algebra has become standard in the literature, and we stick to it.
There is no conflict with the earlier notation $\Delta$ for the diagonal
map of a coalgebra since diagonal maps will not appear explicitly
any more.

It is clear that a generator
determines the  Gerstenhaber bracket.
An observation due to {\it Koszul\/} \cite\koszulon\ 
(p. 261)
says that, for any  Batalin-Vilkovisky algebra
$(\Cal A, [\cdot,\cdot], \partia)$,
the operator $\partia$ 
(which is exact by assumption)
{\it behaves as a derivation for
the  Gerstenhaber bracket\/} $[\cdot,\cdot]$,
that is,
$$
\partia [x,y] = [\partia x,y] -(-1)^{|x|} [x,\partia y], 
\quad x,y \in \Cal A.
\tag3.5.2
$$
A generator $\Delta$, even if exact,
behaves
as a derivation for the multiplication
of $\Cal A$ only if the Gerstenhaber bracket is zero.

\smallskip\noindent
{\smc Definition 3.6.}
Let $(\Cal A,\Delta)$
be a weak  Batalin-Vilkovisky algebra, 
write
$[\cdot,\cdot]$
for the  Gerstenhaber bracket generated by $\Delta$,
and let $d$ be a differential of degree $1$
which endows $(\Cal A,[\cdot,\cdot])$ with a differential  
Gerstenhaber
algebra
structure. 
The triple $(\Cal A,\Delta,d)$
is called a {\it weak differential\/}
Batalin-Vilkovisky algebra provided the graded commutator
${
[d,\Delta] = d\Delta + \Delta d
}$ 
on $\Cal A$ (which is a degree zero operator) is zero.
In particular,
a weak differential
Batalin-Vilkovisky algebra 
$(\Cal A,\partia,d)$
which has $\partia$ exact
is called a 
{\it differential Batalin-Vilkovisky algebra\/}.

\proclaim {Proposition 3.7}
For any 
weak differential  Batalin-Vilkovisky algebra
$(\Cal A,\Delta,d)$,
the differential
$d$
behaves as a derivation for
the  Gerstenhaber
bracket $[\cdot,\cdot]$
on $\Cal A$ generated by $\Delta$, that is to say,
$$
d[x,y] = [dx,y] -(-1)^{|x|} [x,dy], \quad x,y \in \Cal A.
$$
In other words,
$(\Cal A ,[\cdot,\cdot], d)$
is a differential  Gerstenhaber algebra.
\endproclaim

\demo{Proof} This 
is well known; the reader is invited to concoct a proof himself.\qed

Notice that,
under the circumstances of
(3.7),
$\Delta$ need not behave as a derivation for the
 Gerstenhaber bracket
unless $\Delta$ is exact.

Barannikov-Kontsevich
\cite\barakont\ and
Manin
\cite\maninbtw\ (III.10.1.1),
\cite\maninfiv\ (6.1.1),
are  concerned not just with
a Batalin-Vilkovisky algebra
$\Cal A$
for which $\roman H(\fra g)$
is abelian, but moreover with one which
satisfies
the formalization of the
ordinary
formality Lemma
for K\"ahler manifolds.
This leads to a situation of exactly the kind
isolated in Theorem 2.9.
We now explain this in our framework.

We will say that
a differential  Batalin-Vilkovisky algebra
$(\Cal A,\partia,d)$
{\it satisfies the statement of the K\"ahlerian formality lemma\/}
provided the maps
$$
\left(\roman{ker}(\partia),d|_{\roman{ker}(\partia)}\right)
@>{\subseteq}>>
(\Cal A,d),
\quad
\left(\roman{ker}(\partia),d|_{\roman{ker}(\partia)}\right)
@>{\roman{proj}}>>
\roman H (\Cal A,\partia) 
$$
are isomorphisms on homology
where 
$\roman H (\Cal A,\partia)$ is endowed with the zero differential.

What we mean
by the K\"ahlerian formality lemma
is more usually referred to as the
\lq\lq $\partial\overline\partial$-Lemma\rq\rq,
cf. \cite\degrmosu,
but the notation $\partial$
conflicts with our use thereof for a perturbation
to 
which we will stick.

\proclaim{Theorem 3.8}
Let $(\Cal A,\partia,d)$
be a differential  Batalin-Vilkovisky algebra
satisfying the statement of the K\"ahlerian formality lemma,
let $\fra g$ be the differential graded Lie algebra
related with the underlying
differential Gerstenhaber$(\Cal A,[\cdot,\cdot],d)$
as in {\rm (3.2)} above, and extend
the projection $\roman{proj}$
to a contraction
$$
\Nsddata {\fra m} {\nabla}{\pi}{\roman H(\fra g)}h
\tag3.8.1
$$
of chain complexes,
where 
$\fra m =\left(\roman{ker}(\partia),d|_{\roman{ker}(\partia)}\right)$
and $\pi = \roman{proj}$.
Then, as a graded Lie algebra, $\roman H( \fra g)$ is abelian,
and the  data determine a solution 
$\tau \in \roman{Hom}(\Sigm^{\roman c}[s\roman H( \fra g)], \fra g)$
of the master equation
$d \tau  =\frac 12 [\tau,\tau]$
in such a way that 
the following hold:
\newline\noindent
{\rm (i)} The values of $\tau$ lie in $\fra m$,
that is, the composite 
$\Delta \circ \tau
\colon 
\Sigm^{\roman c}[s\roman H( \fra g)] \to \fra g
$
is zero;
\newline\noindent
{\rm (ii)}
the composite $\pi \tau$ coincides with the universal twisting cochain
$\Sigm^{\roman c}[s\roman H( \fra g)] \to  \roman H( \fra g)$
for the abelian graded Lie algebra $\roman H( \fra g)$;
whence
\newline\noindent
{\rm (iii)} for $k \geq 2$, the values of 
the component $\tau_k$ 
on 
$\Sigm_k^{\roman c}[s\roman H( \fra g)]$
of $\tau$
lie in $\roman {im} \Delta$.
\endproclaim

\demo{Proof}
This in an immediate consequence of Theorem 2.9
except statement (iii) which,
in view of (i) and (ii),
follows from the exactness 
of the sequence
$$
0
@>>>
\roman {im} \Delta
@>>>
\roman {ker} \Delta
@>>>
\roman H \fra g
@>>>
0
$$
and the fact that, for
$k \geq 2$,
the component 
on
$\Sigm_k^{\roman c}[s\roman H( \fra g)]$
of the universal twisting cochain
for $\roman H( \fra g)$ is zero, cf. (1.5.5). \qed
\enddemo

Under the circumstances of {\rm (3.8)} suppose that,
in degree 0, 
$\Cal A$ consists of a single copy of the ground ring $R$
necessarily generated by 
the unit 1 of $\Cal A$
and that
$\Delta (1)$ is zero.
Then 1 generates a central copy of $R$ in $\fra g$,
and we may write $\fra g$ as a direct sum
$R \oplus \widetilde{\fra g}$ of differential graded Lie algebras
where
$\widetilde{\fra g}$
is the uniquely determined complement of $R$ in $\fra g$
and $\widetilde{\fra g}$
is itself a differential graded Lie algebra.

\proclaim{Addendum 3.8.2}
Under these circumstances, if the class $[1]$ in homology is non-zero,
the requisite contraction
{\rm (3.8.1)}
may be chosen in such a way that,
for $k \geq 2$,
the values 
of the $\tau_k$
lie in
$\widetilde{\fra g}$.
\endproclaim

\demo{Proof}
The differential graded Lie algebra $\fra m$
decomposes accordingly
as a direct sum
$R \oplus \widetilde{\fra m}$ of differential graded Lie algebras
and so does
the homology of $\fra g$, i.~e. it decomposes as
$
\roman H(\fra g) = R \oplus \roman H(\widetilde{\fra g}).
$
The projection
$\roman{proj}$
from $\roman {ker} \Delta$ to
$\roman H(\Cal A, \Delta)$
may thus be extended
to a contraction of chain complexes
${
\Nsddata {R \oplus \widetilde{\fra m}} {\nabla}{\pi}
{R \oplus \roman H(\widetilde{\fra g})}h
}$
in such a way that the morphisms
$\nabla$, $\pi$and $h$ decompose accordingly;
in particular,
we have a contraction
${
\Nsddata {\widetilde{\fra m}} {\widetilde{\nabla}}{\widetilde{\pi}}
{\roman H(\widetilde{\fra g})}{\widetilde h}.
}$
Let $\widetilde \tau
\colon
\Sigm^{\roman c}[s\roman H( \widetilde{\fra g})] @>>> \widetilde {\fra g}
$
be the twisting cochain
resulting from applying Theorem 2.9 to
this  contraction 
and define
$$
\tau
\colon
\Sigm^{\roman c}[s\roman H(\fra g)]
=
\Sigm^{\roman c}[s\roman H(\widetilde{\fra g})]
\otimes
\Sigm^{\roman c}[sR]
@>>> R \oplus \roman H(\widetilde{\fra g}) = \fra g 
$$
by $\tau = \varepsilon \otimes \widetilde \tau + \tau_0\otimes \varepsilon$
where
$\tau_0$ refers to the universal twisting cochain
$
\tau_0\colon \Sigm^{\roman c}[sR] @>>> R
$
for the abelian Lie algebra $R$
and $\varepsilon$ to the corresponding counits.
The twisting cochain $\tau$ has the desired properties. \qed
\enddemo

Apart from the context,
Theorem 3.8 and the Addendum 3.8.2 contain Lemma 6.1 in \cite\barakont.
More precisely, condition (6.1 (1)) in \cite\barakont,
referred to there as the {\it universality\/} condition, 
translates to the fact that
a contraction is a very precise way of spelling out
a homology isomorphism.
Furthermore, the condition (6.1 (2)) in \cite\barakont,
referred to by the wording {\it flat coordinates\/},
is implied by (i) and (iii) in Theorem 3.8.
The statement of the Addendum amounts to
Barannikov and Kontsevich's 
\lq\lq flat identity\rq\rq\ property.
A statement of the kind (ii) in Theorem 3.8
is referred to by the wording
{\it normalized\/} in Theorem III.9.2 of \cite\maninbtw\ 
and Theorem 4.2 of \cite\maninfiv.
We note that
a contraction  may in fact carry additional information
which could be physically relevant,
e.~g. in terms of higher order correlation functions.

\medskip\noindent{\bf 4. Deformation theory}\smallskip\noindent
Barannikov and Kontsevich set up their construction
in the context of deformation theory.
In particular, they show that one very specific
differential graded Lie algebra
which arises from a Calabi-Yau manifold
is formal
and that 
its {\it formal\/} moduli space 
which, in a sense,
is the {\it formal\/} space
of equivalence classes of
solutions of the relevant master equation,
is affine with a very simple (coordinate) description.
Our method applies much more generally.
Given a differential graded Lie algebra $\fra g$
over the reals or over any field $\fiel k$
of characteristic zero,
the construction of the solution
$\tau$  
given above as (2.8.2) 
can be carried out 
once the requisite contraction 
of the kind (2.2)
has been chosen
(over a field of characteristic zero,
such a choice is always possible, as we remarked earlier)
and still yields a formal
solution of the master equation, but with a perturbed
differential on $\Sigm^{\roman c}[s \roman H \fra g]$,
starting with the operator
induced by the Lie bracket of
$\roman H \fra g$.
The moduli space interpretation is then available
along the lines of that of Schlessinger-Stasheff
\cite\schlstas\
used for the moduli space of rational homotopy types,
once a cohomology algebra is fixed in advance.
One major difference is that the resulting
\lq\lq versal moduli\rq\rq\ space
need no longer be affine.
We now explain this briefly.

Let $L$ be a differential graded Lie algebra 
over a field $\fiel k$ of characteristic zero
which,
viewed as graded by cohomology degrees, is concentrated in nonnegative degrees.
Thus we have homogeneous components
$L^0, L^1, \dots$, and the differential $d$ is an operator
of the form $d\colon L^j \to L^{j+1}$.
The construction in Theorem 2.7 then yields
a 
coalgebra perturbation
$\Cal D$
on
$\Sigm^{\roman c}[S^{-1}\roman H(L)]$
turning the latter into a differential graded coalgebra,
together with a Lie algebra twisting cochain
$$
\tau \colon
\Sigm^{\roman c}_{\Cal D}[S^{-1}\roman H(L)] @>>> L.
$$
Under the circumstances explained in
Theorem 3.8,
$\Cal D$ is zero,
and this reflects the formality of the corresponding
differential graded Lie algebra.
Pursuing the philosophy in \cite\schlstas\
we now offer an interpretation for a non-zero
$\Cal D$.

Recall that a possible deformation theory interpretation
proceeds as follows:
Let $V_L\subseteq L^1$  be the space of 
\lq\lq integrable elements\rq\rq,
or \lq\lq perturbations\rq\rq, or
solutions
$\gamma \in L^1$ of the master equation
$
d \gamma = \frac 12 [\gamma,\gamma].
$
Under suitable circumstances,
$V_L$ is a quadratic variety
(in the appropriate sense).
When 
$\fiel k$ is the field 
of real numbers $\Bobb R$
or that
of complex numbers
$\Bobb C$,
we may consider the  Lie group $\Cal L=\roman{exp}(L^0)$
which corresponds to $L^0$ 
and,
when we assume that the 
induced $\Cal L$-action
on $V_L$
is complete,
we may consider the corresponding
{\it moduli space\/}
$M_L$,
the space of
$\Cal L$-orbits in $V_L$.
Taking an appropriate algebra
$A(L^1)$ of functions 
on the affine space $L^1$ as coordinate ring,
for example polynomial functions with respect to a basis,
or formal power series,
on $V_L$  we have the algebra
$A(V_L) = A(L^1)\big / J$
of functions
where $J$ is the ideal of functions
in
$A(L^1)$ 
which vanish
on $V_L$;
likewise,
the algebra
${
A(M_L) = (A(L^1))^{\Cal L} \big / J^{\Cal L}
}$
of
${\Cal L}$-invariant functions
in $A(L^1)$ modulo the ideal 
$J^{\Cal L}$
of
${\Cal L}$-invariant functions
which vanish
on $V_L$
appears as an algebra
$A(M_L)$ of functions
on the moduli space $M_L$,
and we may view the algebra $A(M_L)$ as a coordinate ring for
$M_L$.
It is {\it not\/} a coordinate ring in the ordinary sense of algebraic
geometry, though.

Since
$\Sigm^{\roman c}_{[\cdot,\cdot]}[S^{-1}(L)]$
is bigraded, so is
$\roman {Hom}(\Sigm^{\roman c}_{[\cdot,\cdot]}[S^{-1}(L)], \Bobb C)$,
and its (co)homology inherits a bigrading.
Consider the bidegree (0,0) cohomology
$$
\roman H^{0,0}
(\roman {Hom}(\Sigm^{\roman c}_{[\cdot,\cdot]}[S^{-1}(L)], \Bobb C)),
$$
that is, the subalgebra of the degree zero cohomology
$\roman H^0(\roman{Hom}
(\Sigm^{\roman c}_{[\cdot,\cdot]}[S^{-1}(L)], \Bobb C))$
which is generated by
cocycles of bidegree (0,0).
The essential observation is now that,
in the formal sense, i.~e. when we are working with
formal power series,
this 
algebra 
is closely related to the algebra
$A(M_L)$
introduced before.
(We leave open here the precise relationship between
$A(M_L)$ and
$\roman H^{0,0}
(\roman {Hom}(\Sigm^{\roman c}_{[\cdot,\cdot]}[S^{-1}(L)], \Bobb C))$.)

Over a general field $\fiel k$ (of characteristic zero),
we now change gears and 
{\it view\/}
the bidegree zero cohomology algebra
$\roman H^{0,0}(\roman{Hom}
(\Sigm^{\roman c}_{[\cdot,\cdot]}[S^{-1}(L)], \fiel k))$
of the (differential graded) $\fiel k$-valued
Chevalley-Eilenberg algebra of $L$
as the coordinate ring
$A(M_L)$ of the moduli space
$M_L$, thereby {\it defining\/}
this moduli space by its coordinate ring.
Theorem 7.1 of \cite\schlstas\ then entails the following interpretation:
By construction,
the morphism
$$
(\Sigm^{\roman c}\nabla)_{\partial}
\colon
\Sigm^{\roman c}_{\Cal D}[S^{-1}(\roman H(L))]
@>>>
\Sigm^{\roman c}_{[\cdot,\cdot]}[S^{-1}(L)]
$$
is 
compatible with the bigradings and
a quasi isomorphism
whence 
$(\Sigm^{\roman c}\nabla)_{\partial}$
induces a cohomology isomorphism
$$
(\Sigm^{\roman c}\nabla)_{\partial}
\colon
\roman H^{*,*}(\roman{Hom}
(\Sigm^{\roman c}_{[\cdot,\cdot]}[S^{-1}(L)], \fiel k))
@>>>
\roman H^{*,*}(\roman{Hom}
(\Sigm^{\roman c}_{\Cal D}[S^{-1}(\roman H(L))], \fiel k)).
$$
Thus our coordinate ring
$A(M_L)$ for $M_L$
now appears as the algebra
\linebreak
$\roman H^{0,0}(\roman{Hom}
(\Sigm^{\roman c}_{\Cal D}[S^{-1}(\roman H(L))], \fiel k))$.
When we reinterpret this in terms of
$\roman H(L)$,
we find that, at least formally,
$M_L$ may be written as the quotient 
${
M_L = W_L\big /F
}$
of what is called the {\it miniversal\/} \lq\lq variety\rq\rq\ 
$W_L \subseteq \roman H^1(L)$
by an equivalence relation $F$.
Here $W_L$ is defined by formal power series
which are determined by the
perturbation
$\Cal D$ and, likewise,
the equivalence relation
$F$ 
is determined by   
$\Cal D$.

In particular,
when $L$ is formal, 
that is, when
$\tau$ and $\Cal D$
can be constructed
in such a way that
$\Cal D$ is merely given by
the induced Lie bracket on $\roman H(L)$,
then
$W_L$ is the pure quadratic \lq\lq variety\rq\rq\ 
which consists of all $\eta \in \roman H^1(L)$
having the property that
$[\eta,\eta] =0 \in \roman H^2(L)$.
Under such circumstances if, in addition,
$\fiel k$ is the field
of real numbers,
the (infinitesimal) $\roman H^0(L)$-action
on $W_L$ then determines a foliation,
and the moduli space
$M_L$ is the space of leaves;
under appropriate circumstances,
$M_L$ may in fact be written as the 
space of $\roman{exp}(\roman H^0(L))$-orbits
with respect to an induced action
of the Lie group
$\roman{exp}(\roman H^0(L))$
on $W_L$.

Under the circumstances of
Barannikov and Kontsevich \cite\barakont,
the induced bracket on $\roman H^*(L)$
is actually trivial whence 
the \lq\lq thickened\rq\rq\ 
moduli space (which initially is a subspace of $\roman H^*(L)$)
is affine, in fact it is all of 
$\roman H^*(L)$;
or for the Kodaira-Spencer algebra
the moduli space would be just
$\roman H^1(L)$.

Schlessinger-Stasheff
\cite\schlstas\
consider deformations of rational homotopy types
which can be described also as deformations
of connected graded commutative algebras
with \lq\lq symmetric\rq\rq\  $A_{\infty}$-structures.
Rather than reviewing that theory, here is a particular application
which also provides a good example of how our result
is more general than that of Barannikov and Kontsevich.

Consider the case of the cohomology of a wedge of spheres
$X= \vee S^{n_j}$. Its cohomology algebra $\roman H^*$ with any field
coefficients $\fiel k$ is 
(in upper degrees) a non-negatively graded vector space
with $\roman H^0 = \fiel k$ and all products  of
positive degree
elements being zero.
The rational homotopy groups 
$\pi_*(\Omega X) \otimes \Bobb Q$
of the based loop space
$\Omega X$
of $X$
are then
isomorphic to
${\Cal L\/}(s\overline {\roman H}_*(X))$
\cite\hiltoone,
the free graded Lie algebra
generated by the reduced homology 
$\overline {\roman H}_*(X)$
of $X$, shifted by 1 in grading.
The deformation 
theory of rational homotopy types refers to the classification 
of rational homotopy types with the same cohomology algebra;
the different types can, in this example,
be distinguished by {\it Massey\/} products: higher order operations
in cohomology.
Attaching a cell by an ordinary Whitehead product
$[S^p,S^q]$
means
the cell carries the product cohomology class.
We rule out this change in the algebra.
Massey and Uehara \cite{\uehamass,\, \massey}
introduced Massey products
in order to detect cells attached by iterated Whitehead products
such as
$[S^p,[S^q,S^r]]$.

Now consider the differential graded Lie algebra
$\roman{Coder}({\Cal L^{\roman c} \/}(S\overline {\roman H}^*(X)))$
where
\linebreak
${\Cal L^{\roman c} \/}(S^{-1}\overline {\roman H}^*(X))$
denotes the free Lie {\it co}algebra on the shifted
cohomology $\overline {\roman H}^*(X)$.
Since all products of positive degree elements are zero,
we start with $d=0$ on
${\Cal L^{\roman c} \/}(S\overline {\roman H}^*(X))$.
A perturbation $\theta$ (of the zero differential)
decomposes into pieces $\theta_k$
which can be identified with
homomorphisms
$$
\theta_k \colon
\roman H^{\otimes (k+2)}
@>>>
\roman H .
$$
Consider the case in which
$\theta = \theta_k$ is homogeneous
and has non-zero image
in only one
$\roman H^n$.
The corresponding space $Y$ is obtained
from $X$ by attaching
the corresponding $n$-cells not as spheres in the bouquet
but non-trivially according to the iterated Whitehead
product determined by the pre-images under $\theta$.
In this space,
all Massey products of order less than $k+2$
will vanish but $\theta$ will represent
a (sum of) non-trivial Massey product(s).

The differential graded Lie algebra
$\pi_*(\Omega Y) \otimes \Bobb Q$
is no longer
sh-equivalent
to the homology of
${\Cal L\/}(s\overline {\roman H}_*(X))$
with trivial differential but rather to the homology with
respect to the corresponding perturbation.
In other words, this differential graded Lie algebra
is not formal.

An example with
\lq\lq continuous moduli\rq\rq,
i.~e., of a one-parameter family of homotopy types,
was first mentioned to us by J. Morgan:
Let $X = S^3 \vee S^3 \vee S^{12}$, so that the space
of possible 5-fold Massey products
$\roman H^{\otimes 5} \to \roman H$
is of dimension 6---equivalently, the attaching maps
are in $\pi_{11}(S^3 \vee S^3) \otimes \Bobb Q$
which is of dimension 6.
A quick check of the relevant dimensions shows that
there are no possible infinitesimal automorphisms,
and 
$\roman{Aut}(\roman H) = \roman{GL}(2) \times \roman{GL}(1)$
is of dimension 5;
thus the Massey products distinguish at least a 1-parameter family.

\bigskip 
\widestnumber\key{999}
\centerline{References}
\smallskip\noindent

\ref \no \barakont
\by S. Barannikov and M. Kontsevich
\paper Frobenius manifolds and formality of Lie algebras of polyvector fields
\paperinfo {\tt alg-geom/9710032}
\jour Internat. Math. Res. Notices
\vol 4
\yr 1998
\pages 201--215
\endref

\ref \no \batviltw
\by I. A. Batalin and G. S. Vilkovisky
\paper Quantization of gauge theories
with linearly dependent generators
\jour  Phys. Rev. 
\vol D 28
\yr 1983
\pages  2567--2582
\endref

\ref \no \ebrown
\by E. Brown
\paper Twisted tensor products.I.
\jour Ann. of Math.
\vol 69
\yr 1959
\pages  223--246
\endref

\ref \no \cartanon
\by H. Cartan
\paper Notions d'alg\`ebre diff\'erentielle; applications aux groupes 
de Lie et aux vari\'et\'es o\`u op\`ere un groupe de Lie
\jour Coll. Topologie Alg\'ebrique
\paperinfo Bruxelles
\yr 1950
\pages  15--28
\endref

\ref \no \cartantw
\by H. Cartan
\paper La transgression dans un groupe de Lie et dans un espace
fibr\'e principal
\jour Coll. Topologie Alg\'ebrique
\paperinfo Bruxelles
\yr 1950
\pages  57--72
\endref

\ref \no \cartanse
\by  H. Cartan
\paper Alg\`ebres d'Eilenberg--Mac Lane et homotopie
\paperinfo expos\'es 2--11
\jour S\'eminaire H. Cartan 1954/55
\publ Ecole Normale Superieure, Paris, 1956
\endref

\ref \no \cartanei
\by H. Cartan and S. Eilenberg
\book Homological Algebra
\publ Princeton University Press
\publaddr Princeton
\yr 1956
\endref

\ref \no \chenone
\by K.T. Chen
\paper Iterated path integrals
\jour Bull. Amer. Math. Soc.
\vol 83
\yr 1977
\pages  831--879
\endref

\ref \no \chenfou
\by K. T. Chen
\paper Extension of $C^{\infty}$
Function Algebra by Integrals and Malcev Completion of $\pi_1$
\jour Advances in Mathematics
\vol 23
\yr 1977
\pages 181--210
\endref

\ref \no \degrmosu
\by P. Deligne, P. Griffiths, J. Morgan, and D. Sullivan
\paper Real homotopy theory of K\"ahler manifolds
\jour Inv. Math.
\vol 29
\yr 1975
\pages  245--274
\endref

\ref \no \eilmactw
\by S. Eilenberg and S. Mac Lane
\paper On the groups ${\roman H(\pi,n)}$. I.
\jour Ann. of Math.
\vol 58
\yr 1953
\pages  55--106
\moreref
\paper II. Methods of computation
\jour Ann. of Math.
\vol 60
\yr 1954
\pages  49--139
\endref

\ref \no \geschthr
\by M. Gerstenhaber and S. D. Schack
\paper Algebras, bialgebras, quantum groups and algebraic
deformations
\paperinfo In: Deformation theory and quantum groups, with
applications to mathematical physics, M. Gerstenhaber and J. Stasheff, eds.
\jour Cont. Math.
\vol 134
\pages 51--92
\publ AMS
\yr 1992
\publaddr Providence 
\endref

\ref \no \getzltwo
\by E. Getzler
\paper Batalin-Vilkovisky algebras and two-dimensional topological field
theories
\jour Comm. in Math. Phys.
\vol 195
\yr 1994
\pages 265--285
\endref

\ref \no \gugenhtw
\by V.K.A.M. Gugenheim
\paper On the chain complex of a fibration
\jour Illinois J. of Mathematics
\vol 16
\yr 1972
\pages 398--414
\endref

\ref \no \gugenhth
\by V.K.A.M. Gugenheim
\paper On a perturbation theory for the homology of the loop space
\jour J. of Pure and Applied Algebra
\vol 25
\yr 1982
\pages 197--205
\endref

\ref \no \gugenlam
\by V.K.A.M. Gugenheim and L. Lambe
\paper Perturbation in differential homological algebra
\jour Illinois J. of Mathematics
\vol 33
\yr 1989
\pages 566--582
\endref

\ref \no \gulasta
\by V.K.A.M. Gugenheim, L. Lambe, and J.D. Stasheff
\paper Algebraic aspects of Chen's twisting cochains
\jour Illinois J. of Math.
\vol 34
\yr 1990
\pages 485--502
\endref

\ref \no \gulstatw
\by V.K.A.M. Gugenheim, L. Lambe, and J.D. Stasheff
\paper Perturbation theory in differential homological algebra. II.
\jour Illinois J. of Math.
\vol 35
\yr 1991
\pages 357--373
\endref

\ref \no \gugenmun
\by V.K.A.M. Gugenheim and H. J. Munkholm
\paper On the extended functoriality of Tor and Cotor
\jour J. of Pure and Applied Algebra
\vol 4
\yr 1974
\pages  9--29
\endref

\ref \no \hiltoone
\by P. J. Hilton
\paper On the homotopy groups of the union of spheres
\jour J. of the London Math. Soc.
\yr 1955
\vol 30
\pages 154--172
\endref

\ref \no \homotype
\by J. Huebschmann
\paper The homotopy type of $F\Psi^q$. The complex and symplectic
cases
\paperinfo
in: Applications of Algebraic $K$-Theory
to Algebraic Geometry and Number Theory, Part II,
Proc. of a conf. at Boulder, Colorado, June 12 -- 18, 1983
\jour Cont. Math.
\vol 55
\yr 1986
\pages 487--518
\endref

\ref \no \perturba
\by J. Huebschmann
\paper Perturbation theory and free resolutions for nilpotent
groups of class 2
\jour J. of Algebra
\yr 1989
\vol 126
\pages 348--399
\endref

\ref \no \cohomolo
\by J. Huebschmann
\paper Cohomology of nilpotent groups of class 2
\jour J. of Algebra
\yr 1989
\vol 126
\pages 400--450
\endref

\ref \no \modpcoho
\by J. Huebschmann
\paper The mod $p$ cohomology rings of metacyclic groups
\jour J. of Pure and Applied Algebra
\vol 60
\yr 1989
\pages 53--105
\endref

\ref \no \intecoho
\by J. Huebschmann
\paper Cohomology of metacyclic groups
\jour Trans. Amer. Math. Soc.
\vol 328
\yr 1991
\pages 1-72
\endref

\ref \no \extensta
\by J. Huebschmann
\paper 
Extensions of Lie-Rinehart algebras and the Chern-Weil construction
\paperinfo in: Festschrift in honor of J. Stasheff's 60th birthday
\jour Cont. Math. 
\vol 227
\yr 1999
\pages 145--176
\publ Amer. Math. Soc.
\publaddr Providence R. I.
\endref

\ref \no \bv
\by J. Huebschmann
\paper Lie-Rinehart algebras, Gerstenhaber algebras, and Batalin-
Vilkovisky algebras
\jour Annales de l'Institut Fourier
\vol 48
\yr 1998
\pages 425--440
\endref

\ref \no \lrbata
\by J. Huebschmann
\paper Rinehart complexes and Batalin-Vilkovisky algebras
\paperinfo preprint 2000
\endref

\ref \no \twilled
\by J. Huebschmann
\paper Twilled Lie-Rinehart algebras and differential Batalin-Vilkovisky 
algebras
\paperinfo math.DG/9811069
\endref

\ref \no \banach
\by J. Huebschmann
\paper Differential Batalin-Vilkovisky algebras arising from
twilled Lie-Rinehart algebras
\paperinfo Poisson Geometry 
\jour Banach Center publications 
\vol 51
\yr 2000
\pages 87--102
\endref

\ref \no \berikas
\by J. Huebschmann
\paper Berikashvili's functor $\Cal D$ and the deformation equation
\paperinfo Fest-\linebreak
schrift in honor of N. Berikashvili's 70th birthday;
{\tt math.AT/9906032}
\jour Proceedings of the A. Razmadze Mathematical Institute
\vol 119
\yr 1999
\pages 59--72
\endref

\ref \no \huebkade
\by J. Huebschmann and T. Kadeishvili
\paper Small models for chain algebras
\jour Math. Z.
\vol 207
\yr 1991
\pages 245--280
\endref

\ref \no \kodnispe
\by M. Kodaira, L. Nirenberg and D. C. Spencer
\paper On the existence of deformations of complex analytic structures
\jour Ann. of Math.
\vol 68
\yr 1958
\pages 450--457
\endref

\ref \no \kosmathr
\by Y. Kosmann-Schwarzbach 
\paper Exact Gerstenhaber algebras and Lie bialgebroids
\jour  Acta Applicandae Mathematicae
\vol 41
\yr 1995
\pages 153--165
\endref

\ref \no \koszulon
\by J. L. Koszul
\paper Crochet de Schouten-Nijenhuis et cohomologie
\jour Ast\'erisque,
\vol hors-s\'erie,
\yr 1985
\pages 251--271
\paperinfo in E. Cartan et les Math\'ematiciens d'aujourd'hui, 
Lyon, 25--29 Juin, 1984
\endref

\ref \no \ladastas
\by T. Lada and J. Stasheff
\paper Introduction to sh Lie algebras for physicists
\jour Int. J. Theor. Phys.
\vol 32
\yr 1993
\pages 1087--1104
\endref

\ref \no \liazutwo
\by B. H. Lian and G. J. Zuckerman
\paper New perspectives on the BRST-algebraic structure
of string theory
\jour Comm. in Math. Phys.
\vol 154
\yr 1993
\pages  613--646
\endref

\ref \no \maninbtw
\by Yu. I. Manin
\book Frobenius manifolds, quantum cohomology, and moduli spaces
\bookinfo Colloquium Publications, vol. 47
\publ Amer. Math. Soc.
\publaddr Providence, Rhode Island
\yr 1999
\endref

\ref \no \maninfiv
\by Yu. I. Manin
\paper Three constructions of Frobenius manifolds: a comparative study
\paperinfo Sir Michael Atiyah: a great mathematician of the
twentieth century,
{\tt Math.QA/9801006}
\jour Asian Math. J.
\vol 3
\yr 1999
\pages 179--220
\endref

\ref \no \massey
\by W. S. Massey
\paper Some higher order cohomology operations
\jour Symposium internacional d\'e topologia algebraica
\paperinfo Universidad Nacional Aut\'onoma de M\'exico and UNESCO, Mexico City
\yr 1958
\pages  145--154
\endref

\ref \no \moorefiv
\by J. C. Moore
\paper Differential homological algebra
\jour Actes, Congr\`es intern. math. Nice
\publ Gauthiers-Villars
\publaddr Paris, 1971
\yr 1970
\pages 335--339
\endref

\ref \no \munkholm
\by H. J. Munkholm
\paper The Eilenberg--Moore spectral sequence and strongly homotopy
multiplicative maps
\jour J. of Pure and Applied Algebra
\vol 9
\yr 1976
\pages  1--50
\endref

\ref \no \quilltwo
\by D. Quillen
\paper Rational homotopy theory
\jour Ann. of Math.
\vol 90
\yr 1969
\pages  205--295
\endref

\ref \no \sanebfou
\by S. Saneblidze
\paper Obstructions to the section problem in fibre bundles
\jour manuscripta math.
\vol 81
\yr 1993
\pages 95--111
\endref

\ref \no \schlstas
\by M. Schlessinger and J. Stasheff
\paper Deformation theory and rational homotopy type
\jour Pub. Math. Sci. IHES
\paperinfo to appear; new version July 13, 1998
\endref

\ref \no \schlsttw
\by M. Schlessinger and J. Stasheff
\paper The Lie algebra structure
of tangent deformation theory 
\jour J. of Pure and Applied Algebra
\vol 38
\yr 1985
\pages 313--322
\endref

\ref \no \uehamass
\by H. Uehara and W. S. Massey
\paper The Jacobi identity for Whitehead products
\jour Algebraic geometry and topology, a symposium in honor of S. Lefschetz
\yr 1957
\pages 361--377
\endref

\ref \no \xuone
\by P. Xu
\paper 
Gerstenhaber algebras and BV-algebras
in Poisson geometry
\jour Comm. in Math. Phys. 
\vol 200
\yr 1999
\pages 545--560
\endref

\end